\input amstex

\def\Aut{\text{Aut}}

\def\b1{\text{\bf 1}}

\def\CA{{\Cal A}}
\def\CB{{\Cal B}}
\def\CC{{\Cal C}}
\def\CD{{\Cal D}}
\def\CF{{\Cal F}}
\def\CE{{\Cal E}}
\def\CH{{\Cal H}}
\def\CI{{\Cal I}}

\def\CM{{\Cal M}}

\def\CO{{\Cal O}}
\def\CP{{\Cal P}}
\def\CR{{\Cal R}}

\def\CV{{\Cal V}}

\def\gr{\text{gr}}

\def\End{\text{End}}
\def\Hom{\text{Hom}}

\def\Ext{\text{Ext}}

\def\#{\,\check{}}

\def\fE{{\frak E}}

\def\id{\text{id}}
\def\tr{\text{tr}}

\def\Coker{\text{Coker}}
\def\Ker{\text{Ker}}

\def\Spec{\text{Spec}}


\def\hra{\hookrightarrow}
\def\iso{\buildrel\sim\over\rightarrow} 

\def\lra{\longrightarrow}


\documentstyle{amsppt}
\NoBlackBoxes

\topmatter\title    Height pairing and nearby cycles       \endtitle \author A.~Beilinson \endauthor   

\dedicatory To   Yuri Ivanovich Manin with deepest gratitude
\enddedicatory

\abstract 
We prove that, as was conjectured by Spencer Bloch, the Hodge period of some limit Hodge structures equals the height pairing of algebraic cycles on the resolution of singularities of the singular  fiber. 
\endabstract

 \subjclass Primary  14C25; Secondary 14D07  \endsubjclass


\address Department of Mathematics, University of Chicago, Chicago, IL 60637\endaddress \email abeilins57\@gmail.com \endemail

\endtopmatter

\document

\head \S1.   Introduction: the theorem and the idea of the proof \endhead

\subhead{ 1.1}\endsubhead {\it The Hodge period.} Suppose we have    a $\Bbb Q$-Hodge structure $E$ with weights in $[-2, 0]$ equiped with isomorphisms  $\iota_0 :  \gr^W_0 E = \Bbb Q (0)$, $\iota_{-2}:  \gr^W_{-2}E = \Bbb Q (1)$.  One defines the   Hodge period $\langle E \rangle =\langle E, \iota_0 ,\iota_{-2} \rangle   \in \Bbb R$  as follows. Consider the $\Bbb R$-Hodge structure $E\otimes \Bbb R$.
 Since the weight filtration on any $\Bbb R$-Hodge structure   with two consequitive weights (canonically)  splits one has $E\otimes \Bbb R = G\oplus \gr^W_{-1} E\otimes \Bbb R$ where $G$ is an extension of $\Bbb R (0)$ by $\Bbb R(1)$. Our $\langle E \rangle$ is the class of this extension in $\Ext^1 (\Bbb R (0),\Bbb R (1))=\Bbb R$.

\remark{Remark}  One computes $\langle E\rangle$ explicitly as follows. Let $E_{\Bbb R}$ be $E\otimes \Bbb R$ viewed a plain $\Bbb R$-vector space,
 $E_{\Bbb C}$ be its complexification. Let $1_{F^0} \in F^0 \subset E_{\Bbb C}$ be any lifting of $\iota^{-1}_0 (1)$. Then $\langle E\rangle$ is the image of $1_{F^0}$ in  $(E_{\Bbb R}+W_{-1}E_{\Bbb C})/(E_{ \Bbb R} +(F^0 \cap W_{-1} E_{\Bbb C})) \buildrel{\sim}\over\leftarrow W_{-2}E_{\Bbb C}/W_{-2}E_{\Bbb R}=\Bbb C/2\pi i\Bbb R \buildrel{\sim}\over\leftarrow \Bbb R$. 
\endremark

\subhead{ 1.2}\endsubhead {\it A geometric example}. Let $Y $ be a smooth proper equidimensional algebraic variety over $\Bbb C$. We denote by $H_i (Y)$ the homology of $H_i (Y(\Bbb C ),\Bbb Q )$ seen as an object of the category of $\Bbb Q$-Hodge structures; ditto for relative homology, etc. Let $Z_m (Y)$ be the group of algebraic $m$-cycles on $Y$ with $\Bbb Q$-coefficients, $Z_m (Y)_0 := \Ker ( cl : Z_m (Y) \to H_{2m}(Y)(-m))$ be the subgroup of cycles homologically equivalent to zero. For a closed subset $P\subset Y$ let $Z_m (P) \subset Z_m (Y)$ be the subgroup of cycles supported on $P$, $Z_m (P)_0 :=Z_m (P)\cap Z_m (Y)_0$. For an $m$-cycle $A$ on $Y$ we denote by $|A|$ its support (which is a closed subset of $Y$).

Suppose $m+m' =\dim Y -1$ and
we have $A\in Z_m (Y)_0$, $B\in Z_{m'}(Y)_0$
such that   $|A|\cap |B|=\emptyset$. Set $E_{|A|,|B|}:= H_{2m+1}(Y\smallsetminus |B|,|A|)(-m)$. Notice that $E_{|B|,|A|}=E_{|A|,|B|}^* (1)$ by the Poicar\'e duality.

\proclaim{ Lemma}  $E_{|A|,|B|}$ has weights in $[-2,0]$. One has $\gr^W_{-2}E_{|A|,|B|} = Z_{m'} (|B| )_0^* (1)$,   $\gr_{-1}^W  E_{|A|,|B|} =H_{2m+1}(Y)(-m)$, $\gr^W_{0}E_{|A|,|B|} = Z_m (|A| )_0$.
\endproclaim
\demo{Proof} Notice that $H_{2m} (|A|)(-m)=Z_m (|A|)$ and $H_{>2m}(|A|)=0$. By  Poincar\'e duality  
 $H_i (Y,Y\smallsetminus |B|)(-\dim Y)=H_{2\dim Y -i}(|B|)^*$, hence $H_{2m+2}(Y,Y\smallsetminus |B|)(-m)= (H_{2m'}(|B|)(-m'))^* (1)  $ and $H_{<2m+2 }(Y,Y\smallsetminus |B|)=0$. Now use the long exact homology sequences for $(Y\smallsetminus |B|,|A|)$ and $(Y,Y\smallsetminus |B| )$. \hfill$\square$
\enddemo

Denote by $E_{A,B}$ the Hodge structure obtained from $E_{|A|,|B|}$ by  pullback by $A$ and pushforward by $B$: 
$$\spreadmatrixlines{1\jot}
\matrix  
Z_n (|B| )_0^* (1)  &  \hra   &  E_{|A|,|B|}
 &\twoheadrightarrow  & Z_m (|A| )_0 
   \\  
B \downarrow &&&&  \uparrow   A
\\
 \Bbb Q (1) &   \hra &      E_{A,B}
 &  \twoheadrightarrow  &  \Bbb Q (0) 
    \endmatrix
\tag 1.2.1$$
 
Our $E_{A,B}$ is as in 1.1, so we have $\langle E_{A,B}\rangle \in \Bbb R$.

\subhead{1.3}\endsubhead {\it The height pairing} (cf.~\cite{B}, \cite{Bl1}).
Let $k$ be a subfield of $\Bbb C$ and suppose  that  $Y$ comes from a variety $Y_k$ over $k$,  $Y=Y_k \otimes \Bbb C$. Let $Z_m (Y_k ) \subset Z_m (Y)$ be the group of algebraic cycles with $\Bbb Q$-coefficients on $Y_k$, $Z_m (Y_k )_0 := Z_m (Y_k )\cap Z_m (Y)_0$, and let $CH_m (Y_k )_0 \subset CH_m (Y_k )$ be their quotients modulo the rational equivalence relation. One checks (see \S 2) that if $A$, $B$ as above are cycles on $Y_k$ then  the class of $\langle E_{A,B}\rangle $ in $ \Bbb R /\Bbb Q \log |k^\times |$ depends only on the  linear equivalence classes of $A$ and $B$, and so one has a bilinear {\it height pairing} $$\langle \,\, ,\,\, \rangle_{Y_k} : CH_m (Y_k )_0 \otimes CH_{m'} (Y_k )_0 \to \Bbb R/\Bbb Q \log |k^\times |. \tag 1.3.1$$
Namely $\langle a,b\rangle_{Y_k}= \langle E_{A,B}\rangle$ where $A$, $B$ are any cycles on $Y_k$ of classes $a$, $b$ such that $ |A|\cap |B|=\emptyset$.

\remark{Remark}  If $k=\Bbb Q$ and we assume some motivic rationality conjectures (see   (2.2.1), (2.2.3) of \cite{B}) then 
$\langle E_{A,B}\rangle$   can be corrected (by adding a finite sum of corrections $\log (p) \langle E_{A,B}\rangle_p$ where $p$ is a prime, $\langle E_{A,B}\rangle_p$  is defined using the Gal$(\Bbb Q_p )$-action on $E_{A,B}\otimes\Bbb Q_\ell$) so that the resulting real number depends only on rational equivalence classes of $A$ and $B$. In this manner (1.3.1) lifts naturally to an $\Bbb R$-valued pairing. 
 \endremark

\subhead{ 1.4}\endsubhead   Finding elements of  Chow groups that are homologically equivalent to zero is an art. Spencer Bloch described one situation where they naturally arise, and conjectured that the height pairing of his cycles can be computed in a different way, namely, as Hodge periods of some nearby cycles. We start with preliminaries.

Let $X$ be a smooth variety over $\Bbb C$ of pure dimension $n\ge 2$, $S$ be a smooth curve, $0\in S$ be a closed point, and $f: X\to S$ be a proper map which is smooth outside a finite subset $\{ x_\alpha \}$ of the fiber $X_0 =f^{-1}(0)$.  Let $Z_\alpha$ be the projectivized tangent cone to $X_0$ at $x_\alpha$; this is a hypersurface in the projectivization $P_{\alpha}:= P (T_{x_\alpha} X )$ of the tangent space; denote by $d_\alpha$ its degree. We assume the next condition:

\enspace 

$(*)$ {\it All hypersurfaces $Z_\alpha$ are smooth.}

\enspace

Let $\pi: Y\to X_0$ be the blowup of $X_0$ at $\{ x_\alpha \}$. Condition $(*)$ implies that $Y$ is a smooth variety, and $Z_\alpha$ are pairwise disjoint divisors on $Y$. Set $Z:= \sqcup Z_\alpha$ and $K:= \Ker (H_{n-2} (Z )\to H_{n-2} (Y))=\text{Im}(H_{n-1}(Y,Z )\to H_{n-2} (Z ))$. If $n=2$ then $Z_\alpha$ is a collection of $d_\alpha$ points; let $K_0  $ be the subgroup of those elements $A=\Sigma A_\alpha \in K\subset \oplus H_0 (Z_\alpha )$ that $\deg A_\alpha =0$ for every $\alpha$.
One has a natural map $H_{n-1}(X_0 )\to H_{n-1}(Y,Z )$ defined as the composition $H_{n-1}(X_0 )\to H_{n-1}(X_0 , \{ x_\alpha \}) \buildrel{\sim}\over\leftarrow H_{n-1}(Y,Z )$.

\proclaim{ Lemma} (i) The map $H_{n-1}(X_0 )\to H_{n-1}(Y,Z )$ is  an isomorphism if $n>2$. If $n=2$ it is injective and its image equals the preimage of $K_0$ in $H_{n-1}(Y,Z )$.
  \newline (ii) $H_{n-1}(Y,Z )$ has weights $1-n$ and $2-n$, and $\gr^W_{2-n}H_{n-1}(Y,Z )=K$. The map $H_{n-1}(Y)\to H_{n-1}(Y,Z )$ has image $W_{1-n}H_{n-1}(Y,Z )$. If $n$ is even then $H_{n-1}(Y)\iso W_{1-n} H_{n-1}(Y,Z )$.
\endproclaim

\demo{Proof} (i) Replace $ H_{n-1}(Y,Z )$ by $H_{n-1}(X_0 , \{ x_\alpha \})$ and use  the long exact homology sequence. (ii) The first assertions follow from the exact homology sequence and purity of weights on $H_\cdot (Y), H_\cdot (Z )$. The last one comes because $H_{n-1}(Z )=0$ if $n$ is even (since $Z_\alpha$ are hypersurfaces).  \hfill$\square$
\enddemo

\subhead 1.5 \endsubhead Consider the variation   of $\Bbb Q$-Hodge structures $\CV$  on $S \smallsetminus \{ 0\}$ with fibers $\CV_s = H_{n-1}(X_s )$. One has a nondegenerate intersection pairing $(\, ,\, ): \CV \otimes\CV \to \Bbb Q (n-1 )$.
Choose a parameter $t$ at $0\in S$ and consider the limiting (a.k.a. nearby cycles) Hodge structure $\psi_t \CV$. Let $\psi^{\text{un}}_t \CV $ be its direct summand where the monodromy acts unipotently. Since  $\psi_t^{\text{un}}$ commutes with duality, $(\, ,\, )$ yields self-duality pairing on it that we denote again by $(\, ,\, )$.  
One has the log of monodromy morphism $N=N_\CV : 
\psi^{\text{un}}_t \CV (1) \to \psi^{\text{un}}_t \CV$ and the specialization morphism $sp:  \psi^{\text{un}}_t \CV \to H_{n-1}(X_0 )$. Let $(\psi^{\text{un}}_t\CV)_N := \Coker (N_\CV )$ be the monodromy coinvariants.  The next assertion follows from the local invariant cycles theorem, see 3.5 for a detailed proof: 

\proclaim{ Proposition} $sp$ factors through the isomorphism $(\psi^{\text{un}}_t \CV)_N \iso H_{n-1}(X_0 )$. 
\endproclaim

\proclaim{Corollary} $\psi^{\text{un}}_t \CV$ has weights in $[-n,2-n]$. One has  $\gr^W_{2-n}\psi^{\text{un}}_t \CV =K$ if $n>2$ and  $\gr^W_{2-n}\psi^{\text{un}}_t \CV =K_0$ if $n=2$. By self-duality, $\gr^W_{-n}\psi^{\text{un}}_t \CV =(\gr^W_{2-n}\psi^{\text{un}}_t \CV )^* (n-1)$. If $n$ is even then $\gr^W_{1-n}\psi^{\text{un}}_t \CV = H_{n-1}(Y)$. 
\endproclaim

\demo{Proof} Since $\psi^{\text{un}}_t \CV$ is self-dual and $N$ is nilpotent, the claim follows from the proposition and the lemma in 1.4. \hfill$\square$ \enddemo

\subhead{ 1.6}\endsubhead {\it Bloch cycles.} We are in the setting of 1.4; suppose $n$ is even, $n=2m+2$. Let $A=\Sigma A_\alpha$ be an $m$-cycle on $Z$. We say that $A$ is a  {\it Bloch cycle} if it is homologically equivalent to zero on $Y$, i.e., $cl (A)$ lies in $K(-m)\subset H_{n-2}(Z )(-m)$. If $m=0$ then we demand, in addition, that $cl (A)\in K_0\subset K$. 

\proclaim{Lemma} If $A$ is a Bloch cycle then each  $cl (A_\alpha )\in H_{n-2}( Z_\alpha )(-m)$ is primitive.
\endproclaim

\demo{Proof} The composition $H_{n-2}(Z )(-m) \to H_{n-2}(Y)(-m) \to H_{n-4}(Z_\alpha )(-m+1)$, where the second arrow is the pullback by $Z_{\alpha}\hra Y$, sends any class $c=\Sigma c_\alpha$ to $c_\alpha \cap c_1 (\CO (-1))$ (for $\CO (-1)$ is the normal bundle to $Z_\alpha$ in $Y$). This composition kills $cl(A_\alpha )$ since the first arrow does. \hfill$\square$
\enddemo

If $A,B$ are two Bloch cycles then we denote by $E^\psi_{A,B}=E^\psi_{A,B,t}$ the Hodge structure obtained from $\psi^{\text{un}}_t \CV (-m)$ by pullback by $cl (A)$ and pushforward by $cl (B)^*$:
$$\spreadmatrixlines{1\jot}
\matrix  
K^* (m+1)  &  \to   &   \psi^{\text{un}}_t \CV (-m)
 &\to  & K(-m)
  \\  
cl(B)^* \downarrow &&&&  \uparrow   cl(A)
\\
\,\,\,\,\,\,\,\,\,\,\, \Bbb Q (1) &  \,\,\, \hra &      E^\psi_{A,B}
 & \!\!\! \twoheadrightarrow  & \!\!\!\!\!\!\!\! \Bbb Q (0) 
    \endmatrix
\tag 1.6.1$$

Our $E^\psi_{A,B} $ is as in 1.1 so we have  $\langle E^\psi_{A,B} \rangle \in \Bbb R$.  

\subhead 1.7 \endsubhead {\it Examples.} Consider the case when  we have a single singular point $x_0 \in X_0 $ of $f$ and the singularity at $x_0$ is quadratic. Then the monodromy action on $\psi_t \CV$ is unipotent,  the only possible Bloch cycle is the difference $A$ of the rulings of the quadric $Z_0$, and it is actually a Bloch cycle if and only if the monodromy action on $\psi_t \CV$ is nontrivial or, equivalently, the Hodge structure on $H_{n-1}(X_0 )$ is not pure. 
\proclaim{Lemma}
 (i) If $m=0$ then the curve $X_0$ can have either 1 or 2 irreducible components, and $A$ is a Bloch cycle if and only if $X_0$ is irreducible. \newline (ii) If $X/S$ is a family of quadratic hypersurfaces in $\Bbb P^{n}$ then $A$ is not a Bloch cycle.  \newline (iii) If $X/S$ is a family of hypersurfaces of degree $d$ on a given smooth projective variety $P$ then $A$ is a Bloch cycle if $d$ is large enough.
\endproclaim

\demo{ Proof} (i) is clear. (ii) follows since the global monodromy for quadratic hypersurfaces is $\pm 1$, and so it can't contain non-trivial unipotent local monodromy. 

(iii) Consider the corresponding map $r: S\to B:=\{$hypersurfaces of degree $d$ on $P \}$.
Since $X$ is smooth $r$ is transversal to the locus $D\subset B$ of degenerate hypersurfaces. Replacing $S$ by a germ of another transversal to $D$ that intersects $D$ near $r(0)$ would not change the topology of $X$ over a small disc around $0$. So we can 
 assume that $S$ is a Zariski open subset of the base of a Lefschetz pencil on $P$. Then, since  local monodromies of a Lefschetz pencil are all conjugate,\footnote{Recall that this follows from  irreducibility of the dual variety of $X$.} triviality of one local monodromy amounts to triviality of the global monodromy. Thus $A$ is a Bloch cycle if and only if the global monodromy on $\CV$ is not trivial. Let us check that this happens for large enough $d$.

 If $R\subset P$ is the axis of our pencil then $H^\cdot (X)= H^\cdot (P) \oplus H^{\cdot-2}(R)(-1)$, and so $h^{n-1,0} (P)= h^{n-1,0} (X)$ which equals $h^{n-1, 0}(X_s )$ if the global monodromy is trivial. Thus the monodromy is not trivial when $h^{n-1, 0}(X_s )>h^{n-1,0} (P)$. To finish the argument it remains to notice that
 $h^{n-1, 0}(X_s ) \ge \dim ( H^0 (P,\Omega_P^{n}(d))/H^0 (P, \Omega_P^{n}))$, and so it tends to $\infty$ when
$d\to \infty$.
 \hfill$\square$ \enddemo

\subhead{ 1.8}\endsubhead {\it Statement of the theorem.}
Now suppose  we have a subfield $k\subset \Bbb C$ and our datum is defined over $k$, i.e., there is $X_k / S_k$, a closed point $0$ of $S_k$, 
a parameter $t$ on $S_k$ at $0$, and Bloch cycles $A$, $B$ on $Z_{k}$ such that $X/S$, etc., come by base change $k\to \Bbb C$. Let $a\in CH_m (Y_k )_0$, $b\in CH_{m}(Y_k )_0$ be the classes of $A$ and $B$.
 The next result was conjectured  by Spencer Bloch:

\proclaim{ Theorem} One has $\langle a,b\rangle_{Y_k} = \langle E^\psi_{A,B} \rangle \! \mod \Bbb Q \log |k^\times |$.
\endproclaim

In case $n=1$ the theorem was proven in \cite{BlJS}.

\remark{Remark} Suppose we are in the situation of Remark in 1.3. If $\langle E^\psi_{A,B}\rangle$ is corrected in the same way as was discussed there, then the theorem lifts to an equality of real numbers. The proof does not change; we will not discuss it below. 
\endremark

\subhead{ 1.9}\endsubhead {\it Reformulation of the theorem that discards Hodge periods; the idea of the proof.} Let $A'$, $B'$ be cycles on $Y_k$ of classes $a$, $b$ such that $|A'|\cap |B'|=\emptyset$ (notice that they are, most probably, not supported on $Z_k$).
We want to show that $\langle E_{A',B'}\rangle =\langle E^\psi_{A,B}\rangle$ (see 1.2, 1.6).
Let us compare the Hodge structures $E=E_{A',B'}$ and $E^\psi =E^\psi_{A,B}$ themselves. Their weights lie in $[-2,0]$, and one has a canonical identification $\gr^W_\cdot E = \gr^W_\cdot E^\psi$. Indeed, 
 $\gr^W_{0}E^{(\psi )}=\Bbb Q (0)$, $\gr^W_{-2}E^{(\psi )}=\Bbb Q (1) $ by  the constructions, and  $\gr^W_{-1} E =H_{2m+1}(Y)(-m)=\gr^W_{-1} E^{(\psi )}$ by the lemma in 1.2, and the one in 1.4 combined with the corollary in 1.5.   This identification lifts (uniquely) to $W_{-1}E = W_{-1}E^\psi$ and $E/W_{-2}E =E^\psi /W_{-2}E^\psi$. Indeed, the classes of
  extensions $0\to H_{2m+1} (Y)(-m) \to E^{(\psi )}/W_{-2} E^{(\psi )}\to \Bbb Q (0)\to 0$ both
 equal Deligne cohomology class $cl_{\CD}(A)$  (a.k.a. Griffiths' Abel-Jacobi periods) of $A$; by duality, the classes of (the duals to) extensions $0\to \Bbb Q (1)\to W_{-1} E^{(\psi )}\to H_{2m+1} (Y)(-m) \to 0$  both equal to $cl_{\CD}(B)$ (see loc.cit.).

Now suppose we have a $\Bbb Q$-Hodge structure $H$ of weight $-1$ and two classes $a\in \Ext^1 (\Bbb Q (0),H)$, $b\in \Ext^1 (H,\Bbb Q (1))$. Consider the set $\fE^{\CH}_{a,b} =\fE^{\CH}(H)_{a,b} $ of all Hodge structures $E$ with weights in $[-2,0]$  and equipped with identifications $\gr^W_0 E =\Bbb Q (0)$, $\gr_{-1}^W E=H$, $\gr^W_{-2}E=\Bbb Q(1)$ such that the extensions $E/W_{-2}E$ and $W_{-1}E$ have classes $a$ and $b$. The group $\Ext^1 (\Bbb Q(0),\Bbb Q (1))=\Bbb C^\times \otimes\Bbb Q$ acts on $\fE^{\CH}_{a,b}$ by the Baer sum action, and $\fE^{\CH}_{a,b}$ is a $\Bbb C^\times \otimes\Bbb Q$-torsor. Notice that for $q\in \Bbb C^\times $ one has $\langle q \cdot E\rangle = \log |q| +\langle E\rangle$. Applying this format to $H= H_{2m+1}(Y)(-m)$, $a=cl_{\CD}(A)$, $b=cl_{\CD}(B)$ and
 $E_{A',B'}, E^\psi_{A,B}\in \fE^{\CH}_{a,b} $ we get $E_{A',B'} - E^\psi_{A,B}\in \Bbb C^\times \otimes \Bbb Q$. Now the theorem in 1.8 follows immediately from the next result (notice that the Hodge periods and the height pairing play no role here):

\proclaim{ Theorem} One has $E_{A',B'} - E^\psi_{A,B}\in k^\times \otimes \Bbb Q \subset \Bbb C^\times \otimes \Bbb Q$.
\endproclaim
 
The theorem would be an immediate corollary of the motivic formalism if all the above constructions could  be spelled in motivic world: Indeed, we would have then a motivic version $\fE^{\CM}$ of $\fE^{\CH}$ which is  an $\Ext^1_{\CM}(\Bbb Q (0),\Bbb Q(1)) =k^\times\otimes \Bbb Q$-torsor equipped with the Hodge realization embedding  $\fE^{\CM}\hra \fE^{\CH}$;  our $E_{A',B'}$, $ E^\psi_{A,B}$ would come from elements of $\fE^{\CM}$, and so their difference lies in $k^\times\otimes \Bbb Q$. The only problem is that in the present day formalism of motives, due to Voevodsky, Ayoub, and Cisinski-D\'eglise, the t-structure is  not available, so we do not have the motivic version of separate homology groups like $H_i (Y)$. The actual proof is an exercise in spelling out the constructions in a way that makes the t-structure  redundant.

\remark{Warning} (I owe it to the referee.) The proof  uses compatibility of the Hodge realization of motivic sheaves with the six functors formalism. It is not fully checked yet, so the theorem is conditional subject to that  verification (see 4.1).
\endremark

\enspace

I am very grateful to Spencer Bloch for explaining me his conjecture and stimulating discussions (pity  Spencer refused to    coauthor   the article), to Volodya Drinfeld for valuable comments and discussions,  to Luc Illusie for calling my attention to the construction of  \cite{I} which helped to clearify and simplify  the argument, and to the referee for his remarks and the above Warning.

\head \S 2. The height pairing and the construction of $E^\CM_{a,b}\in \fE_{a,b}^{\CM}\subset \fE^{\CH}_{a,b}$
\endhead

This section is a variation on the theme of \cite{Bl2} and \cite{G}.

\subhead{2.1} \endsubhead
Let $\CC$ be a stable dg category. 
It yields two other dg categories $\CC^{(1)}$ and $\CC^{(2)}$ constructed as follows: 

An object of  $\CC^{(1)}$ is a closed morphism $\alpha : M\to N$ of degree 0 in $\CC$. One has $\Hom ((M,N,\alpha ), (M',N',\alpha' ))^i = \Hom (M,M')^i \times \Hom (N,N' )^i \times \Hom (M,N')^{i+1}\subset \Hom (\CC one (\alpha ),\CC one (\alpha' ))^i$, and the differential is defined so that the latter embedding is a morphism of complexes; the composition of morphisms is defined in a similar way. There are three dg functors
$\CC^{(1)}\to \CC$ which send $(M,N,\alpha )$ to $M$, $N$, and $\CC one (\alpha )$ respectively. We can view $\CC^{(1)}$ as the category of distinguished triangles, and the rotation yields an autoequivalence $\rho :\CC^{(1)}\to \CC^{(1)}$ which sends $\alpha :M \to N$ to $\rho (\alpha ): N \to \CC one (\alpha )$; the inverse autoequivalence is $\rho^{-1}(\alpha ): \CC one (\alpha )[-1]\to M$.

An object of $\CC^{(2)}$ is a datum $(P,M,Q, \alpha ,\beta ,\kappa )$ where $P$, $M$, $Q$ are objects of $\CC$,  $\alpha \in \Hom (P, M )^1$,   $\beta\in \Hom (M, Q)^1 $  are closed maps, and $\kappa\in \Hom (P,Q )^1$ is such that $d(\kappa )=\beta \alpha$; we sometimes abbreviate it to $(\alpha,\beta,\kappa )$. One can assign to such a datum an object $E= E(\alpha,\beta,\kappa  )\in \CC$ which equals $P\oplus M\oplus Q$ with $\alpha $, $\beta $, and $-\kappa$ added as the components to the differential.\footnote{Thus  $E=\CC one ((\alpha,\kappa):P[-1]\to \CC one (\beta : M[-1]\to Q ))=\CC one ((\kappa,\beta):\CC one (\alpha :P[-2]\to M[-1])\to Q )$.} There is a filtration $Q\subset \CC one (\beta : M[-1],Q) \subset E$, and morphisms in 
$\CC^{(2)}$ are the same as morphisms between the corresponding objects $E$ that preserve this filtration.
We have 
 two dg functors $\CC^{(2)}\to \CC^{(1)}$ which send to $(\alpha,\beta,\kappa )$ to $\alpha : P[-1]\to M$ and $\beta : M\to Q[1]$, and six dg functors $\CC^{(2)}\to \CC$ which send $(\alpha,\beta,\kappa )$ to $P$, $M$, $Q$, $\CC one (\alpha :P[-1]\to M )$, $\CC one (\beta : M[-1] \to Q)$, and $E(\alpha,\beta,\kappa )$ respectively. 

The dg category $\CC^{(2)}$ carries a natural involution $\sigma $ which sends $   (P,M,Q,\alpha,\beta,\kappa )$ to the object $ (Q[-1], E(\alpha,\beta,\kappa ), P[1], \alpha^\sigma ,\beta^\sigma  , 0)$ where
$\alpha^\sigma $ and $\beta^\sigma$ are the evident embedding and projection.

\remark{Remark} One can view an object $(\alpha,\beta,\kappa )\in \CC^{(2)}$ as an object of $\CC$ equipped with a 3-step filtration in two different ways. Namely, this could be $E(\alpha ,\beta ,\kappa  )$ equipped with an evident  filtration with successive quotients $Q$, $M$, and $P$. Or this could be $M$ equipped with a filtration whose successive quotients are $P[-1]$, $E(\alpha ,\beta ,\kappa )$, and $Q[1]$. The involution $\sigma$ exchanges the two perspectives. 
\endremark

\subhead{2.2} \endsubhead  For $\CC$ as above we denote by $\CC^\times$ the $\infty$-groupoid of its homotopy equivalences,  by $\CC^{\times \tau}$ the corresponding 1-truncaded plain groupoid, and by $H\CC$ the homotopy category of $\CC$. For
 $S,T\in \CC$ set $\Ext^i (S,T):=H^i \Hom (S,T)=\Hom_{H\CC}(S,T[i])$.  
Denote by $\CE xt (S,T)$ the plain Picard groupoid of extensions that corresponds to the two-term complex $\tau^{[0,1]}\Hom (S,T)$.

For $M,N\in \CC$ let $\CC^{(1)\times}_{M,N}$ be the $\infty$-groupoid of collections $(\alpha' :M'\to N' , \iota_M ,\iota_N )$ where $(\alpha' : M'\to N' )\in \CC^{(1)}$ and $\iota_M :M\to M'$, $\iota_N :N\to N'$ are homotopy equivalences. It is equivalent to  the Picard $\infty$-groupoid that corresponds to the complex $\tau^{\le 0}\Hom (M,N)$. The 1-truncated plain Picard groupoid $\CC^{(1)\times\tau}_{M,N}$  corresponds to the two-term complex $\tau^{[-1,0]}\Hom (M,N)$. 

Similarly, 
for three objects  $P,M,Q\in \CC$ we have  the $\infty$-groupoid $\CC^{(2)\times}_{P,M,Q}$ whose objects are data $(P',M',Q',\alpha',\beta',\kappa', \iota_P ,\iota_M ,\iota_Q )$ where $
(P',M',Q',\alpha',\beta',\kappa')\in \CC^{(2)}$ and $\iota_P :P\to P'$, $\iota_M :M\to M'$, $\iota_Q :Q\to Q'$ are homotopy equivalences. 
 The 1-truncated plain groupoid
$\CC^{(2)\times\tau}_{P,M,Q}$  contains a normal subgroup $\Ext^0 (P,Q)=\Hom_{H\CC}(P,Q )$. Let
  $\fE =\fE (M) =\fE (P,M,Q)$ be the  quotient groupoid. It is equivalent to the groupoid of triples $ (\alpha ,\beta, \kappa )$ where $\alpha \in \Hom (P, M )^1$,   $\beta\in \Hom (M, Q)^1 $  are closed maps, and $\kappa\in \Hom (P,Q )^1/d(\Hom (P,Q )^0 )$ is such that $d(\kappa )=\beta \alpha$; a morphism $ (\alpha ,\beta, \kappa ) \to (\alpha' ,\beta', \kappa' ) $ in $\fE$ is a pair $(\phi ,\psi )$ where $\phi \in \Hom (P,M )^0 /d(\Hom (P,M )^{-1})$, $\psi\in \Hom (M,Q )^0 /d(\Hom (M,Q )^{-1})$ are such that $\alpha' -\alpha = d(\phi )$, $\beta' -\beta = d(\psi )$, $\kappa'-\kappa =  \beta \phi  +\psi \alpha +\psi d(\phi )$.

The projection 
$\CC^{(2)}_{P,M,Q}\to \CC^{(1)}_{P[-1], M}\times \CC^{(1)}_{M,Q[1]}$ yields a map of plain groupoids $
\fE (P,M,Q)\to\CC^{(1)\times\tau}_{P[-1], M}\times \CC^{(1)\times\tau}_{M,Q[1]}=\CE xt (P,M)\times\CE xt (M,Q )$, $(\alpha ,\beta,\kappa )\mapsto (\alpha ,\beta )$. The group $\Ext^1 (P,Q)$ acts   on $\fE$ by translations of $\kappa$, and non-empty fibers  $\fE_{\alpha ,\beta }$ are   $\Ext^1 (P,Q)$-torsors. 

\remark{Remark} $\fE (P,M,Q)$ is naturally functorial with respect to $P$ and $Q$: every pair of closed morphisms $\mu: P_1 \to P$ and $\nu: Q\to Q_1$ yields a map $\fE (P,M,Q)\to \fE (P_1 ,M,Q_1 )$, $(\alpha ,\beta,\kappa )\mapsto (\alpha\mu,\nu\beta,\nu\kappa\mu )$; is compatible with the  $\Ext^1 (P,Q)$-action via the map $(\mu^* ,\nu_* ): \Ext^1 (P,Q)\to \Ext^1 (P_1 ,Q_1 )$.

Suppose $\Ext^2 (P,Q)=0$. Then $\fE_{\alpha ,\beta }$ are non-empty, and 
the addition maps $\fE_{\alpha_1,\beta}\times \fE_{\alpha_2,\beta}\to\fE_{\alpha_1 +\alpha_2 ,\beta}$, $\fE_{\alpha,\beta_1 }\times \fE_{\alpha,\beta_2}\to \fE_{\alpha ,\beta_1 + \beta_2}$ define on $\fE$  the structure of an  $\Ext^1 (P,Q)$-biextension of $(\CE xt (P,M),\CE xt (M,Q))$.

\subhead {2.3} \endsubhead In our first example $\CC $ is the dg category whose homotopy category is the bounded derived category $D\CH$ of the category $\CH$ of $\Bbb Q$-Hodge structures, and $P=\Bbb Q(0)$, $Q=\Bbb Q(1)$. We denote the corresponding $\fE$ by $\fE^\CH =\fE^\CH (M)$. Then $\Ext^{\neq 1}_{D\CH} (P,Q)=0$ and $\Ext^1_{D\CH} (P,Q)=\Bbb C^\times \otimes \Bbb Q$, so
 $\fE^\CH$ is a $\Bbb C^\times \otimes \Bbb Q$-biextension of   $(\CE xt (\Bbb Q(0 ),M), \CE xt (M,\Bbb Q (1)))$.

 Let $\Ext^1_0 (\Bbb Q (0),M)\subset \Ext^1 (\Bbb Q (0),M)$, $\Ext^1_0 (M, \Bbb Q (1))\subset  \Ext^1 (M, \Bbb Q (1)) $  be the subgroups of those elements $a$, $b$ such that the maps $H^0 a  : \Bbb Q (0)\to H^1 M$, $H^{-1}b : H^{-1}M \to \Bbb Q (1)$  vanish. Let $\CE xt_0  (\Bbb Q (0),M) \subset \CE xt  (\Bbb Q (0),M)$, etc., be the Picard groupoids of such extensions.

\proclaim  {Lemma} Suppose that $\Hom (\Bbb Q (0), H^0 M)=\Hom (H^0 M,\Bbb Q(1)) =0$. 
\newline
(i) The restriction of $\fE^\CH$ to $(\CE xt_0  (\Bbb Q (0),M), \CE xt_0 (M, \Bbb Q (1)))$ descends to the $\Bbb C^\times \otimes \Bbb Q$-biextension of   $(\Ext^1_0 (\Bbb Q(0 ),M), \Ext^1_0 (M,\Bbb Q (1)))$.
\newline  
 (ii) $\fE^{\CH  }$ is naturally functorial with respect to $M$: if $\varphi : M\to M'$ is a morphism, and we have $a' \in \Ext^1_0 (\Bbb Q (0),M')$, $b'\in \Ext^1_0 (M',\Bbb Q (1))$ with $\varphi_* (a )= a'$, $\varphi^* (b' )=b$ then there is a canonical identification  $\fE^{\CH  }(M)_{a,b}=\fE^{\CH  }(M')_{a',b'}$.
 \newline
(iii) The isomorphisms $\Ext^1_0 (\Bbb Q (0),M)\iso \Ext^1 (\Bbb Q (0),H^0 M)$,
$\Ext^1_0 (M, \Bbb Q (1))\iso 
\Ext^1 $ $(H^0 M, \Bbb Q (1))$ which assign to an extension its zero cohomology,
   lifts naturally to an isomorphism of biextensions $H^0 : \fE^{\CH } (M)\iso \fE^{\CH} (H^0 M)$. One has $EH^0 = H^0 E$.
\endproclaim

\demo{Proof} Let us prove (i); the rest is clear. We need to check that for every closed  $\alpha \in \Hom^1_0 ( \Bbb Q (0),M)$, $\beta \in \Hom^1_0 (M,\Bbb Q(1))$ the action of $\Aut (\alpha  )\times \Aut(\beta )= \Hom ( \Bbb Q (0),M)$ $ \times  \Hom (M,\Bbb Q(1))$ on $\fE^\CH_{\alpha,\beta}$ is trivial.

Since $\CH$ has homological dimension 1 our $M$ is isomorphic to the direct sum of its homologies and so $ \Aut (\alpha )=\Ext^1 (\Bbb Q(0), H^{-1}M)$, $\Aut (\beta ) = \Ext^1 (H^1 (M),\Bbb Q(1))$ by the condition on $M$. The action of $(e,h)\in \Ext^1 (\Bbb Q(0), H^{-1}M)\times \Ext^1 (H^1 (M),\Bbb Q(1))$ on $\fE^\CH_{\alpha,\beta}$ is the translation by $H^{-1}(\beta )e+ hH^0 (\alpha )$ which is 0  since $\alpha ,\beta \in \Ext^1_0$.
\hfill$\square$
\enddemo

Consider the  $\Bbb R$-biextension   $\log |\fE^{\CH } (M) |$ which is the pushforward of  the $\Bbb C^\times \otimes\Bbb Q$-biextension $\fE^{\CH } (M)$ along the homomorphism $\log |\cdot |: \Bbb C^\times \to \Bbb R$.

 \proclaim {2.4. Lemma} Suppose that $H^0 M$ is pure of weight $-1$ (which implies the condition of  the lemma in 2.3). Then the function $\fE^\CH (M)\to \Bbb R$, $(\alpha ,\beta ,\kappa )\mapsto \langle E (\alpha ,\beta ,\kappa ) \rangle := \langle H^0 E (\alpha ,\beta ,\kappa )\rangle$, see 1.1, is a natural trivialization of $\log |\fE^{\CH } (M) |$. \endproclaim

\demo { Proof} Everything said in   2.3 works for the category $\CH_{\Bbb R}$ of $\Bbb R$-Hodge structures. 
The extension of scalars functor $\CH \to \CH_{\Bbb R}$, $?\mapsto ?\otimes \Bbb R$, yields a morphism of our biextensions $\fE^{\CH  } (M)\to \fE^{\CH_{\Bbb R}} (M\otimes \Bbb R )$. The map $\Ext^1 (\Bbb Q (0),\Bbb Q(1)) \to
\Ext^1 (\Bbb R (0),\Bbb R (1))$ equals  $\log |\, |$ after the standard identifications of the $\Ext$ groups with, respectively, $\Bbb C^\times\otimes\Bbb Q$ and $\Bbb R$. Since 
$\Ext^1 (\Bbb R (0),H^0 M\otimes \Bbb R )=\Ext^1 (H^0 M\otimes\Bbb R, \Bbb R(1))=0$ by the condition on $M$, one has $\fE^{\CH_{\Bbb R}} (M\otimes \Bbb R )=\fE^{\CH_{\Bbb R}} (H^0 M\otimes \Bbb R )=\Bbb R$. The map  $\fE^{\CH  } (M)\to \fE^{\CH_{\Bbb R}} (M\otimes \Bbb R )=\Bbb R$ is $\langle \, \rangle$ of 1.1. \hfill$\square$
\enddemo

\subhead 2.5 \endsubhead  Let $k\subset \Bbb C$ be a subfield. Denote by $D\CM (k)$ the dg category of geometric Voevodsky $\Bbb Q$-motives over $k$. We have the Hodge realization dg functor $D\CM (k)\to D\CH  $, $M\mapsto M^{\CH}$. Consider the story of 2.2 for $\CC =D\CM (k)$ with $P=\Bbb Q(0)$, $Q=\Bbb Q(1)$. As before
 one has $\Ext^{\neq 1}_{D\CM (k)} (\Bbb Q (0),\Bbb Q(1))=0$, and there is a canonical identification
$\Ext^1 (\Bbb Q (0),\Bbb Q(1))=k^\times\otimes\Bbb Q$ such that the Hodge realization map between the $\Ext^1$'s is the embedding $k^\times \otimes \Bbb Q \hra \Bbb C^\times \otimes \Bbb Q$. 
So for any $M\in D\CM (k)$ we get a   $ k^\times\otimes\Bbb Q$-biextension of $(\CE xt^1 (\Bbb Q (0),M),\CE xt^1 (M,\Bbb Q (1)))$ together with the Hodge realization morphism $
\fE^{\CM}(M)  \to \fE^{\CH }(M) := \fE^{\CH }(M^{\CH  }) $ of the biextensions.

\remark {Remark} Since the   homomorphism  $k^\times \otimes \Bbb Q \hra \Bbb C^\times \otimes \Bbb Q$    is injective, the maps of torsors $\fE^{\CM}(M)_{\alpha ,\beta} \to \fE^{\CH  }(M)_{\alpha ,\beta}:=  \fE^{\CH  }(M)_{\alpha^\CH ,\beta^\CH }   $ are injective too. \endremark
\enspace

We define $\Ext^1_0 (\Bbb Q (0),M) \subset \Ext^1_0 (\Bbb Q (0),M)$ and $\Ext^1_0 (M, \Bbb Q (1))\subset  \Ext^1 (M, \Bbb Q (1)) $ as preimages of the $\Ext^1_0$ subgroups of the Hodge setting by the Hodge realization maps. Assume that $H^0 M^{\CH}$ is pure of weight $-1$. Then (i) and (ii) of the lemma in 2.3 remain true in the $D\CM (k)$ setting (with $\Bbb C^\times$ replaced by $k^\times$): this follows from loc.cit.~by Remark above. Thus we have a $k^\times \otimes \Bbb Q$-biextension $\fE^{\CM} (M)$ of $( \Ext^1_0 (\Bbb Q (0),M),  \Ext^1_0 (M, \Bbb Q (1)))$ together with a map of biextensions $\fE^{\CM}(M) \to \fE^{\CH} (M)$, so the lemma in 2.4 provides a natural trivialization of the $\Bbb R$-biextension $\log |\fE^{\CM}(M) |$. The image of $\fE^\CM_{a,b}$ in $\Bbb R/\Bbb Q\log |k^\times |$   depends only on $a,b \in \Ext^1_0 (M, \Bbb Q (1))\times \Ext^1_0 (\Bbb Q (0),M) $, and we denote it by $ \langle a,b \rangle_M$. It is clearly biadditive with respect to $a,b$.\footnote{Indeed, a morphism from a biextension by a trivial group to a trivialized biextension amounts to a biadditive pairing.}  We have defined
 a canonical {\it height pairing}    $$\langle \, \,\rangle_M  :\Ext^1_0 (\Bbb Q (0),M) \times \Ext^1_0 (M, \Bbb Q (1))   \to \Bbb R/\Bbb Q\log |k^\times | . \tag 2.5.1$$

\subhead 2.6 \endsubhead   We return to the situation of 1.3 and set $M:=M(Y_k )(-m) [-1-2m ]$ where $M(Y_k )$ is the motive of $Y_k$.
 One has $\Ext^1 (\Bbb Q (0), M)=CH_m (Y_k )$,  $\Ext^1 (M,\Bbb Q(1))=CH_{m'} (Y_k ) $ by  Poincar\'e duality,  and $\Ext^1_0$ are the subgroups $CH (Y_k )_0$ of cycles homologically equivalent to zero. Therefore we get a $k^\times \otimes\Bbb Q$-biextension $\fE^{\CM}$ of $(CH_m (Y_k )_0 , CH_{m'} (Y_k )_0 )$, the map of biextensions $\fE^{\CM}\to \fE^{\CH}$,  
the trivialization of $\log |\fE^\CM |$, and the height pairing $\langle\, ,\,\rangle_M  :
CH_m (Y_k )_0 \times CH_{m'}(Y_k )_0 \to \Bbb R / \Bbb Q\log |k^\times |$.   

By (iii) of the lemma in 2.3 one has $H^0 : \fE^{\CH} (M)\iso  \fE^{\CH} (H_{2m+1}(Y)(-m))$. For $a\in CH_m (Y_k )_0$, $b\in CH_{m'}(Y_k )_0$ pick, as in 1.3,  cycles $A$, $B$ that represent them such that
 $|A|\cap |B|=\emptyset$.\footnote{Recall that $|A|, |B|\subset Y_k$ are supports of the cycles.} Let us construct $(a, b, \kappa_{A,B} ) \in \fE^{\CM}_{a,b}$ such that the Hodge realization 
$E^\CH_{A,B}$ of
  $E^\CM_{A,B}:= E (a,b,\kappa_{A,B})$ (see 2.1) has  zero cohomology $H^0 E^\CH_{A,B}$ equal to the Hodge structure  
 $E_{A,B}$ from 1.3. This would imply that for our $M$ the height pairing (2.5.1) equals 
(1.3.1).

The composition of the maps $M(|A|)\buildrel\alpha\over\to M(Y_k ) \buildrel\beta\over\to M(Y_k ,Y_k \smallsetminus |B|)$ is naturally homotopic to 0: indeed, $M(Y_k ,Y_k \smallsetminus |B|):= \CC one (M(Y_k \smallsetminus |B|)\to M(Y_k))$, and the homotopy $\kappa_{|A|,|B|}$ is $M(|A|)\to M(Y_k \smallsetminus |B|) \subset \CC one$. Thus we have $(  \alpha ,\beta,\kappa_{|A|,|B|})\in D\CM^{(2)}$ (see 2.1). Notice that $E (\alpha ,\beta ,\kappa_{|A|,|B|}) =M(Y_k \smallsetminus |B|,|A|)$.

One has $\Ext^{-2m}(\Bbb Q (m), M(|A|)) =Z_m (|A|):=$ the group of $m$-cycles supported on $|A|$ (recall that $\dim |A|=m$), and $\Ext^{2m+2} (  
M(Y_k ,Y_k \smallsetminus |B|), \Bbb Q (m+1))=Z_{m'}(|B|)$ by  Poincar\'e duality. Therefore we have $(\alpha A, B\beta ,B\kappa_{|A|,|B|} A)=(\Bbb Q (m)[2m +1] ,M(Y_k ),   \Bbb Q (m)[2m +2],    \alpha A, B\beta ,B\kappa_{|A|,|B|} A)     \in D\CM^{(2)}$. The promised $(a, b, \kappa_{A,B} ) \in \fE^{\CM}_{a,b}$ is $(\alpha A, B\beta ,B\kappa_{|A|,|B|} A)(-m)[-1-2m]$. The fact that
$H^0 E^\CH_{A,B}$ equals the Hodge structure  
 $E_{A,B}$ from 1.3 follows from the construction.

\head \S 3. The unipotent nearby cycles in the Hodge setting
\endhead

\subhead 3.1\endsubhead {\it A nearby cycles reminder.}   In this section we play with algebraic varieties over $ \Bbb C$.
For an algebraic variety $X$ we denote by  $\CH (X)$ the abelian category of perverse Hodge $\Bbb Q$-sheaves of M.~Saito on $X $, by $D\CH (X)$ its bounded derived category. It satisfies the usual Grothendieck six functors formalism.  Below ${}^*$ is  Verdier duality. Every object of $\CH (X)$, hence of $D\CH (X)$, carries a canonical weight filtration.

For $\CF \in D\CH (X)$ let $\Gamma (X,\CF),\Gamma_c (X,\CF)    \in D\CH$ be the complex of chains, resp. chains with compact support, with coefficients in $\CF$ equipped with the natural Hodge structure,  $H^\cdot_{(c)} (X,\CF ):= H^\cdot \Gamma_{(c)} (X,\CF )\in \CH$; 
set $\Gamma_{(c)} (X):=\Gamma_{(c)} (X ,\Bbb Q (0)_X ) $, $H^\cdot_{(c)} (X):= H^\cdot_{(c)}  (X,\Bbb Q (0) )$, and denote by $(\, ,\, )$   Poincar\'e duality pairing. Similarly for a closed subvariety $A\subset X$ we set $\Gamma_A (X):=\Gamma_{A} (X,\Bbb Q (0) )\in D\CH $, etc.

Let  
$g: X\to \Bbb A^1$ be a function on $X$; set $X_0 := g^{-1}(0)$, and let $v: X\smallsetminus X_0 \hra X$, $i_{X_0}:X_0 \hra X$ be the open and closed embeddings.
One has the unipotent nearby cycles 
functor $ \psi^{\text{un}}_{g} : D\CH (X\smallsetminus X_0 )\to D\CH (X_0 )$ that carries a natural logarithm of monodromy morphism $N=N_g =N_{\CF}: \psi_g^{\text{un}} (\CF )(1) \to \psi_g^{\text{un}} (\CF )$ where $\CF\in D\CH (X\smallsetminus X_0 )$. It has \' etale local origin with respect to $X_0$.  For sheaves on $X$ there is a natural morphism of functors $\iota: i^*_{X_0} \to  \psi^{\text{un}}_{g} v^*$.
There are basic canonical identifications:
\newline
(i)  Compatibility with Verdier duality: One has $\psi_g^{\text{un}}(\CF^* )= (\psi_g^{\text{un}}\CF )^* (1)[2]$.
\newline
(ii) Compatibility with proper direct images: Suppose  $h: X\to T$ is a proper map and $t$ is a function on  $T$ such that $g=th$; then one has $\psi_t^{\text{un}} h_* \CF =h_{*} \psi^{\text{un}}_{g}  \CF$.
\newline
(iii) One has $\CC one (N_{\CF})= i_{X_0}^* v_* \CF (1)[1]$.
\newline
(iv) For every $n>0$ one has   $\psi^{\text{un}}_{g^n} \CF \iso \psi^{\text{un}}_{g} \CF$.

These identifications are mutually compatible; (i) and (ii) are compatible with the action of $N$, and (iv) identifies $N_{g^n}$ with $nN_g$. Finally, one has 
\newline
(v) $\psi^{\text{un}}[-1]$ is t-exact for the perverse t-structure.

\remark{Examples} Suppose that $X$ is smooth of dimension $n$ and $\CF =\Bbb Q (0)_{X\smallsetminus X_0}$. Then $\CF^* =\CF (n)[2n]$ hence 
$\psi_g^{\text{un}}(\CF )^* = (\psi_g^{\text{un}}\CF ) (n-1)[2n-2]$.
\newline
(a) If $g$ is smooth then $\iota_{\Bbb Q (0)_X} : \Bbb Q (0)_{X_0}\iso \psi^{\text{un}}_g \CF  $, $N_\CF =0$.
\newline
(b) Suppose $g$ is semi-stable and $X_0$ has two irreducible components $Y$ and $Y'$. By (a) 
one has natural morphisms $j_{Y'\smallsetminus Y!} \Bbb Q_{Y'\smallsetminus Y} \to \psi^{\text{un}}_g \CF \to j_{ Y\smallsetminus Y'  *} \Bbb Q_{Y\smallsetminus Y'}$   compatible with the $N$-action (we take it that on the left and right object $N$ acts trivially). They form an exact triangle; its Verdier dual is the same triangle with $Y$ and $Y'$ interchanged.
\endremark

\subhead 3.2\endsubhead We are in the setting of 1.4 and follow the notation there. 

Let $j: U:= X_0 \smallsetminus \{ x_\alpha\} \hra X_0 \hookleftarrow  \{ x_\alpha \} :\sqcup i_{x_\alpha}   $ be the complementary open and closed embeddings.   Let $\CI$ be the intersection cohomology sheaf $j_{!*}\Bbb Q (0)_U =\tau^{\le n-2} j_* \Bbb Q (0)_U$\footnote{Below $\tau$ is the usual truncation, ${}^p \tau$ is the perverse one.} on $X_0$; set $\CI^+ :=\pi_* \Bbb Q (0)_Y$. One has natural  self-duality isomorphisms  $\CI^*=\CI (n-1)[2n-2]$, $\CI^{+*}=\CI^+ (n-1)[2n-2]$ (recall that $Y$ is smooth of dimension $n-1$ and $\pi$ is proper).

The decomposition theorem for $\pi$ is easy and explicit:

\proclaim{Proposition} There is a natural orthogonal direct sum decomposition $$\CI^+ =\CI  \oplus\, \oplus_\alpha i_{x_\alpha  *} \tau^{[2, 2n-4]}\Gamma (P_{\alpha})  \tag 3.2.1$$ compatible with the self-dualities.
\endproclaim
\demo{Proof} One has a natural orthogonal direct sum decomposition $$\Gamma (Z_\alpha )= H^{n-2}_{\text{prim}}(Z_\alpha )[2-n] \oplus \tau^{\le 2n-4} \Gamma (P_{\alpha}) \tag 3.2.2$$ defined as follows. Consider the embedding $Z_\alpha \hra P_\alpha$. The  pullback and Gysin maps $\Gamma (P_\alpha )\to \Gamma (Z_\alpha )\to \Gamma (P_{\alpha} )(1)[2]$ are mutually dual for  Poincar\'e duality pairings, and their composition in either direction equals to the multiplication by $c_1 (\CO (d_\alpha ))$.\footnote{Since $\CO (d_\alpha )$ is the normal bundle to $Z_\alpha$ in $P_\alpha$.} Thus the composition of  $\tau^{\le 2n-4} \Gamma (P_\alpha ) \to \Gamma (Z_\alpha )\to \tau^{\ge 0}( \Gamma (P_\alpha )(1)[2])$ is an isomorphism. This yields a direct sum decomposition $\Gamma (Z_\alpha )$ $= ? \oplus \tau^{\le 2n-4} \Gamma (P_\alpha )$. Since multiplication by  $c_1 (\CO (d_\alpha ))$ preserves the direct sum decomposition,  the only nonzero cohomology of ? is $  H^{n-2}_{\text{prim}}(Z_\alpha )\subset H^{n-2}(Z_\alpha )$, q.e.d.

  Consider the embeddings of smooth divisors  $i_{Z_\alpha} : Z_\alpha \hra Y$. One has $i^!_{Z_\alpha} \Bbb Q (0)_Y =\Bbb Q(-1)[-2]_{Z_\alpha}$, $i^*_{Z_\alpha} \Bbb Q (0)_Y =\Bbb Q(0)_{Z_\alpha}$, and the composition of the adjunction maps $i_{Z_\alpha *}i^!_{Z_\alpha}\Bbb Q (0)_Y \to \Bbb Q (0)_Y \to i_{Z_\alpha *}i^*_{Z_\alpha}\Bbb Q (0)_Y$ equals the multiplication by $c_1 (\CO (-1))$ map $\Bbb Q (-1)[-2]_{Z_\alpha}\to \Bbb Q (0)_{Z_\alpha}$.\footnote{Since $\CO (-1 )$ is the normal bundle to $Z_\alpha$ in $Y$.} Apply $\pi_*$; then  $i^!_{x_\alpha}\CI^+ = \Gamma (Z_\alpha )(-1)[-2]$, $i^*_{x_\alpha}\CI^+ = \Gamma (Z_\alpha )$ by base change, and the composition of the adjunctions $ i_{x_\alpha *} i^!_{x_\alpha}\CI^+ \to\CI^+ \to i_{x_\alpha *} i^*_{x_\alpha}\CI^+$ is multiplication by $c_1 ( \CO (-1))$ map $i_{x_\alpha *}\Gamma (Z_\alpha )(-1)[-2]\to i_{x_\alpha *}\Gamma (Z_\alpha )$. 

Composing  the maps
 $\tau^{\le 2n-6}\Gamma (P_\alpha )  \hra \Gamma (Z_\alpha )  $ and   $\Gamma (Z_\alpha )\twoheadrightarrow \tau^{[2, 2n-4]} \Gamma (P_\alpha ) $ that come from decomposition (3.2.2)   from the left and from the right with the latter adjunctions, we get the maps $i_{x_\alpha *} (\tau^{\le 2n-6}\Gamma (P_\alpha ))(-1)[-2] \to \CI^+ \to i_{x_\alpha *}  \tau^{[2, 2n-4]} \Gamma (P_\alpha ) $. Their composition is an isomorphism, which yields a decomposition $\CI^+ =\CI^? \oplus i_{x_\alpha *}  \tau^{[2, 2n-4]} \Gamma (P_\alpha )  $. Since the adjunctions are mutually dual, the decomposition is orthogonal. 

 By (3.2.2) one has $i^!_{x_\alpha} \CI^? = H^{n-2}_{\text{prim}}(Z_\alpha )(-1)[-n]\oplus \Bbb Q (1-n)[2-2n] ,$ $ i^*_{x_\alpha} \CI^? =  H^{n-2}_{\text{prim}}(Z_\alpha )[2-n] \oplus \Bbb Q (0)$.

Thus $\CI^?  [n-1]$ is a perverse sheaf which equals $\Bbb Q(0)[n-1]_U$ on $U$ and has no subquotients supported on $\{ x_\alpha\}$, and so $\CI^? =\CI$.  We are done. \hfill$\square$
\enddemo

\remark {Remarks} (i) The adjunction map $\Bbb Q (0)_{X_0}\to \pi_* \Bbb Q (0)_Y =\CI^+$ takes value in  $\CI \subset \CI^+$ since $\Hom (\Bbb Q (0)_{X_0}, i_{x_\alpha  *} \tau^{[2, 2n-4]}\Gamma (P_{\alpha}))=0$. \newline (ii) Set $\CB :=  \oplus i_{x_\alpha  *} H^{n-2}_{\text{prim}}(Z_\alpha )[1-n]$. By the formula for $i^*_{x_\alpha} \CI$ at the end of the previous paragraph, one has  an exact triangle $\Bbb Q (0)_{X_0}\to \CI \to \CB [1]$.
\endremark

\subhead {3.3}\endsubhead 
 As in 1.5, $t$ is a local coordinate at $0\in S$; shrinking $S$ we can assume that $t$ is defined and invertible on $S\smallsetminus \{ 0\}$, so $X_0 = (tf)^{-1}(0)$. Consider the functor  $ \psi^{\text{un}}_{tf}: D\CH (X\smallsetminus X_0 )\to D\CH (X_0 )$   (see 3.1). Set  $\CR := \psi^{\text{un}}_{tf} \Bbb Q(0)_{X\smallsetminus X_0 } $. By 3.1(i) one has a canonical self-duality identification $\CR^* =\CR (n-1)[2n-2]$ and the  mutually dual maps $\Bbb Q (0)_{X_0}\buildrel{\iota}\over\to \CR 
\buildrel{\iota^*}\over\to \Bbb Q (0)^*_{X_0}(1-n)[2-2n]$ which are isomorphisms over $U$.

The next result is due to  Illusie \cite{Il}; we will need it in 4.5. The reader can skip it at the moment and jump directly to section 3.4.

\proclaim{Proposition} For every critical point $x_\alpha$ one has canonical isomorphisms   $$i^!_{x_\alpha}\CR = \Gamma_c (P_\alpha \smallsetminus Z_\alpha ), \quad i^*_{x_\alpha}\CR = \Gamma (P_\alpha \smallsetminus Z_\alpha ) \tag 3.3.1$$ interchanged by  duality. The $N$-action on $i^!_{x_\alpha}\CR, i^*_{x_\alpha}\CR$ is trivial.
\endproclaim

\demo{Proof} (a) The claim is local at $x_\alpha$, so for the proof we remove from $X$ the rest
of critical points, and still call it $X$ by the abuse of notation. Let $S^\flat \to S$ be the covering of degree $d_\alpha$ obtained by adding $t^\flat =t^{1/d_\alpha}$ to the sheaf of functions; its Galois group is $\mu_{d_\alpha}$. Set $X^\flat := X\times_S S^\flat$ and let $f^\flat : X^\flat \to S^\flat$ be the projection. Our  $X^\flat $  is a hypersurface  
$\{ (x,t^\flat ): (tf)(x)- t^{\flat d_\alpha} =0\}$ in $X\times \Bbb A^1$; its only singular point   is $(x_\alpha ,0)$. 
The projectivized tangent cone $Q_\alpha $ of $X^\flat$ at $(x_\alpha ,0 )$ is a hypersurface in $
P_\alpha^+ := P(T_{(x_\alpha ,0 )} X\times \Bbb A^1 )$. The Galois group $\mu_{d_\alpha}$ acts on $X^\flat$ hence on $Q_\alpha$.

(b) Let us check that $Q_\alpha$ is a $\mu_{d_\alpha}$-covering of $P_\alpha$ completely ramified along $Z_\alpha$ and \'etale over its complement, and  $Q_\alpha$ is smooth. 
To see this, consider the leading term
 $[tf]_{d_\alpha} (x)$  (of the Taylor expansion) of $tf$ at $x_\alpha$; then the  leading term  of  $(tf)(x)- t^{\flat d_\alpha}$ at $(x_\alpha ,0 )$ is $[tf]_{d_\alpha} (x) - t^{\flat d_\alpha}$. The zeros of $[tf]_{d_\alpha}$ is  $Z_\alpha \subset P_\alpha$, of $[tf]_{d_\alpha}   (x)- t^{\flat d_\alpha}$ is $Q_\alpha \subset P_\alpha^+$, and so the 
 projection  $Q_\alpha \to P_\alpha$ $(x,t^\flat )\mapsto x$, is as claimed. The smoothness of $Q_\alpha$ follows from that of $Z_\alpha$.

(c) Let $\pi^+ : X^+ \to X^\flat$ be the blowup of $X^\flat$ at $(x_\alpha ,0)$. By (b) $X^+$ is smooth and the map $f^+ :=f^\flat \pi^+ : X^+ \to S^\flat$ has semistable reduction at $0\in S^\flat$. The fiber $X^+_0$ has two irreducible components:  one equals $Y$ and the other $Q_\alpha$, and their intersection equals $Z_\alpha$. The action of $\mu_{d_\alpha}$ on $X^\flat$ yields one on $X^+$. The $\mu_{d_\alpha}$-action on $X^+_0$ fixes $Y$ and acts on $Q_\alpha$ as described in (b). The projection $\pi^+_0 : X^+_0 \to X^\flat_0 = X_0$ contracts $Q_\alpha$ to $x_\alpha$. 

Set $\CR^+ := \psi^{\text{un}}_{tf^+} \Bbb Q (0)_{X^+ \smallsetminus X^+_0}$, $\CR^\flat := \psi^{\text{un}}_{tf^\flat} \Bbb Q (0)_{X^\flat \smallsetminus X^\flat_0}$. These are  sheaves on $X^+_0$ and $X^\flat_0 =X_0$ respectively that are  naturally  $\mu_{d_\alpha}$-equivariant. By 3.1(ii) (with $h=\pi^+$) one has a natural identification $\pi^+_{0*} \CR^+ =  \CR^\flat$ compatible with the  $\mu_{d_\alpha}$-actions. Since the projection $p: X^\flat \to X$ is a $\mu_{d_\alpha}$-torsor over $X\smallsetminus X_0$ one has $\Bbb Q (0)_{X \smallsetminus X_0}=(p_* \Bbb Q (0)_{X^\flat \smallsetminus X^\flat_0})^{\mu_{d_\alpha}}$, and so, by
3.1(ii) with $h=p$, one has $\CR = \CR^{\flat \mu_{d_\alpha}}$. Therefore $\CR = (\pi^+_{0*} \CR^+ )^{\mu_{d_\alpha}}. $

(d) By 3.1(iv) with $g=t^\flat f^+$, $n=d_\alpha$, one has $\psi^{\text{un}}_{tf^+} = \psi^{\text{un}}_{t^\flat f^+}$. Our $t^\flat f^+$ is semi-stable, so we have the exact 
triangle $  j_{Y\smallsetminus Z_\alpha  ! }\Bbb Q_{Y\smallsetminus Z_\alpha} \to \CR^+  \to   j_{Q_\alpha \smallsetminus Z_\alpha * }\Bbb Q_{Q_\alpha \smallsetminus Z_\alpha}  $   as in Example (b) in 3.1. Applying  $\pi^+_{0*}$   we get an exact triangle $j_!  \Bbb Q_{U} \to \CR^\flat  \to i_{x_\alpha *} \Gamma (Q_\alpha \smallsetminus Z_\alpha )$. Passing to $\mu_{d_\alpha}$-invariants we get, by (b), an exact triangle $j_!  \Bbb Q_{U } \to \CR  \to i_{x_\alpha *} \Gamma (P_\alpha \smallsetminus Z_\alpha )$; here we use the identification $\Gamma (Q_\alpha \smallsetminus Z_\alpha )^{\mu_{d_\alpha}} \iso
 \Gamma  (P_\alpha \smallsetminus Z_\alpha )$ defined as the composition $\Gamma (Q_\alpha \smallsetminus Z_\alpha )^{\mu_{d_\alpha}}   \subset \Gamma (Q_\alpha \smallsetminus Z_\alpha )\buildrel{\tr}\over\to  \Gamma (P_\alpha \smallsetminus Z_\alpha ) $.  Thus we get the   isomorphism $ i^*_{x_\alpha }\CR \iso \Gamma (P_\alpha \smallsetminus Z_\alpha ) $ in (3.3.1). The second isomorphism there comes in the dual manner from the dual exact triangle $j_{Q_\alpha \smallsetminus Z_\alpha ! }\Bbb Q_{Q_\alpha \smallsetminus Z_\alpha} \to \CR^+  \to  j_{Y\smallsetminus Z_\alpha  * }\Bbb Q_{Y\smallsetminus Z_\alpha} $. Since $\pi^+_{0*}$ commutes with duality, the two isomorphisms are mutually dual, and we are done. \hfill$\square$
\enddemo

Let $\alpha_\CR$ be the composition $ \CB \buildrel{\partial}\over\to \Bbb Q (0)_{X_0}\buildrel{\iota}\over\to \CR $ where $\partial$ is the boundary map of the triangle from Remark (ii) in 3.2, so $\CI =\CC one (\partial )$. Let us compute the map $i^!_{x_\alpha} (\alpha_\CR )$.   Consider the standard triangle $H^{n-2}_{\text{prim}}(Z_\alpha )[1-n] \buildrel{\delta}\over\to \Gamma_c (P_\alpha \smallsetminus Z_\alpha ) \buildrel{\tr}\over\to \Bbb Q(1-n)[2-2n]$ that comes from (3.2.2).  

\proclaim{Lemma} $-i^!_{x_\alpha} (\alpha_\CR )$ equals the composition $\delta_\CR$ of the maps $H^{n-2}_{\text{prim}}(Z_\alpha )[1-n] \buildrel{\delta }\over\to \Gamma_c (P_\alpha \smallsetminus Z_\alpha ) \buildrel{(3.3.1)}\over{=} i^!_{x_\alpha}\CR   $.
\endproclaim

\demo{Proof}  Consider the exact triangle $$j_{Q_\alpha \smallsetminus Z_\alpha !}\Bbb Q(0)_{Q_\alpha \smallsetminus Z_\alpha }\oplus j_{Y \smallsetminus Z_\alpha !}\Bbb Q(0)_{Y\smallsetminus Z_\alpha }\to \Bbb Q(0)_{X_0^+} \to \Bbb Q(0)_{Z_\alpha}. \tag 3.3.2$$ Let $(\delta_Q ,\delta_Y ): \Bbb Q(0)_{Z_\alpha}[-1]\to j_{Q_\alpha \smallsetminus Z_\alpha !}\Bbb Q(0)_{Q_\alpha \smallsetminus Z_\alpha }\oplus j_{Y \smallsetminus Z_\alpha !}\Bbb Q(0)_{Y\smallsetminus Z_\alpha }$ be the boundary map. Its composition with the map to $\Bbb Q(0)_{X_0^+}$, and hence with the further composition with
 $ \Bbb Q(0)_{X_0^+} \buildrel{\iota}\over\to \CR^+$, is 0.   
Therefore the sum of the compositions   
$\Bbb Q(0)_{Z_\alpha}[-1] \buildrel{\delta_Q}\over\lra 
j_{Q_\alpha \smallsetminus Z_\alpha !}  \buildrel{\iota}\over\to \CR^+$ and 
$\Bbb Q(0)_{Z_\alpha}[-1] \buildrel{\delta_Y}\over\lra 
j_{Y \smallsetminus Z_\alpha !}  \buildrel{\iota}\over\to \CR^+$ is 0. Apply $i^!_{x_\alpha}\pi_*^+$ and consider the restriction of our compositions to $ H^{n-2}_{\text{prim}}(Z_\alpha )[1-n] \subset \Gamma (Z_\alpha )[-1]$. For the first one it is $\delta_\CR$, for the second one it is $ i^!_{x_\alpha}( \alpha_\CR )$, and we are done. \hfill$\square$
\enddemo

\subhead 3.4
\endsubhead   Set $\CP :=\CR [n-1]=\psi^{\text{un}}_{tf} \Bbb Q(0)_{X\smallsetminus X_0 }[n-1] $; this  is   a perverse sheaf   on $X_0$; one has a canonical self-duality identification $\CP^* =\CP (n-1)$. Consider the perverse sheaves $\CP^N := \Ker (N: \CP \to \CP (-1))$, $\CP_N := \Coker (N: \CP  (1)\to \CP )$.

\proclaim{Lemma} (i) $\Bbb Q (0)_{X_0}[n-1]$ is a perverse sheaf of weights $n-1$ and $n-2$ with $\gr^W_{n-1} =\CI [n-1]$, $\gr^W_{n-2} =\oplus_\alpha \, i_{x_\alpha *} H^{n-2}_{\text{prim}}(Z_\alpha )$. \newline
(ii) One has $\CP^N = \Bbb Q (0)_{X_0}[n-1]$, $\CP_N = (\Bbb Q (0)_{X_0}[n-1])^* (1-n)$.
\newline (iii) $\CP$ has weights in $[n-2,n]$. One has $W_{n-1}\CP =\Bbb Q (0)_{X_0}[n-1] $, $\CP /W_{n-2}\CP = (\Bbb Q (0)_{X_0}[n-1])^* (1-n) $, $\gr^W_{n-2}\CP = \oplus_\alpha \, i_{x_\alpha *} H^{n-2}_{\text{prim}}(Z_\alpha )$, $\gr^W_{n-1}\CP=\CI [n-1]$, $\gr^W_{n}\CP =   (\gr^W_{n-2}\CP )^* (1-n)$.    
\endproclaim

\demo{Proof} (i) The exact triangle from Remark (ii) in 3.2 amounts to an exact triangle  $\oplus i_{x_\alpha  *} H^{n-2}_{\text{prim}}(Z_\alpha ) \to   \Bbb Q (0)_{X_0} [n-1] \to \CI [n-1]$, and we are done since its left and right terms are pure perverse sheaves of weights $n-2$ and $n-1$ respectively.

 (ii) For any sheaf $\CA$ on $X$ one has a canonical exact triangle $i_{X_0}^* \CA \to i_{X_0}^* v_* v^*\CA \to i^!_{X_0}\CA[1]$: Indeed,  the map $v_! v^*\CA \to v_*  v^* \CA$ factors as composition $v_! v^*\CA \to\CA\to v_*  v^* \CA$, and so one has an exact triangle $\CC one (v_! v^*\CA \to\CA )\to \CC one (v_! v^*\CA \to v_*  v^* \CA )\to \CC one (\CA\to v_*  v^* \CA)$ which is supported on $X_0$. The promised exact triangle is its restriction to $X_0$. 

Now take for $\CA $ the perverse sheaf $\Bbb Q(0)_{X}  [n]$. The first term of the triangle is $\Bbb Q(0)_{X_0}  [n]$ which is perverse sheaf shifted by 1, its third term is $(\Bbb Q (0)_{X_0}[n-1])^* (-n)$ which is a perverse sheaf. Therefore they equal, respectively, ${}^p H^{-1}$ and ${}^p H^0$ of  $i_{X_0}^* v_* v^*\Bbb Q (0)_X[2n] $, i.e., of $ \CC one (N :\CP \to \CP (-1))$ by   3.1(iii), and we are done.

 (iii) Since $N$ is nilpotent, the weights of $\CP$ are bounded from below by the minimum of weights of $\CP^N$, which is $n-2$ by (ii) and (i). By self-duality of $\CP$ they are bounded then from above by $n$, and we have the first assertion. It implies that $W_{n-2}\CP \subset \CP^N$. The rest follows directly from (i), (ii), and self-duality of $\CP$. 
 \hfill$\square$
\enddemo

\subhead {3.5}\endsubhead {\it Proof of the proposition in 1.5.} We use the notation in loc.cit. Injectivity of $sp: (\psi^{\text{un}}_t \CH)_N \to H_{n-1}(X_0 )$ follows from the local invariant cycles theorem. Let us check the surjectivity. By 3.1(ii) applied to $h=f$
 (recall that $f $ is proper) and 3.1(v) applied to $\psi^{\text{un}}_t$, one has $\psi^{\text{un}}_t \CH = H^{0}  (X_0 ,\CP )(n-1)$. By 3.4 we have exact sequence of perverse sheaves $0\to \oplus_\alpha \, i_{x_\alpha *} H^{n-2}_{\text{prim}}(Z_\alpha )(n-1)\to \CP (n-1) \to (\Bbb Q (0)_{X_0}[n-1])^* \to 0$. Its left term has finite support, and so has no  cohomology in degrees $\neq 0$. Therefore  the map $   H^0  (X_0 ,\CP ) (n-1) \to H^0 (X_0 ,  (\Bbb Q (0)_{X_0}[n-1])^* )=   H_{n-1}(X_0 )$ is surjective. This map equals $sp$, and we are done. \hfill$\square$

\head \S4. The motivic setting and the construction of $E^{\psi \CM}_{a,b}\in \fE_{a,b}^{\CM}$
\endhead

\subhead 4.1 \endsubhead  We are in the setting of 1.8 so $k\subset \Bbb C$ is a subfield and we play with varieties over $k$. Changing slightly the notation of 1.3 and 1.8, for a variety $Z=Z_k$ we set
 $Z_{\Bbb C}:= Z \otimes_k \Bbb C$.
The notation of \S3 is preserved except that {\it we equip from now on all Hodge sheaves and Hodge structures  met previously with extra upper index ${}^\CH$.}

We play with  motives (a.k.a.~motivic sheaves) over varieties, see \cite{A1} and \cite{CD}. For a variety $Z$ the category of constructible $\Bbb Q$-motives over $Z$ is denoted by $D\CM (Z )$. We use Grothendieck's six functors formalism for $D\CM$ which was developed, following the ideas of Voevodsky,  in \cite{A1}   and further treated in \cite{CD}. Recall that $D\CM (\Spec\, k )$$=D\CM (k)$ is the category of Voevodsky's geometric $\Bbb Q$-motives over $k$. For  a variety $Z$ one has $M(Z )=\pi_{Z !}\pi_{Z}^! \Bbb Q(0)$ where $\pi_{Z} : Z \to \Spec\, k$ is the structure map. For a motivic sheaf $\CF$ on $Z$ set $\Gamma (Z ,\CF ):= \pi_{Z *} \CF ,  \Gamma_c (Z ,\CF ):= \pi_{Z !} \CF  \in D\CM (k)$; we write $\Gamma  (Z ):=\Gamma  (Z ,\Bbb Q (0)_{Z} )$, ditto for $\Gamma_c (Z)$. For a smooth $Z$ of dimension $d$ one has $\pi_{Z}^! \Bbb Q (0)=\Bbb Q (d)_{Z}[2d]$, and so  $M(Z )=\Gamma_c (Z ) (d)[2d]$.

There is a Hodge realization functor $D\CM (Z ) \to D\CH  (Z_{\Bbb C} )$, $\CF \mapsto \CF^\CH$, that was defined in \cite{Iv}. {\it Below I will assume its
compatibility with the six functors; this was not fully checked yet (though there is little doubt that the claim is true) which makes the story below conditional subject to verification of the compatibility.}

The formalism of unipoteny nearby cycles in the setting of motivic sheaves was developed in \S\S 3.4, 3.6 of  \cite{A2}. The motivic version of everything said in 3.1 holds except property (v) (for the t-structure is not available). The Hodge realization functor commutes with the nearby cycles functors.

\subhead 4.2 \endsubhead {\it Notation:} Notice that $\Hom (\Bbb Q(i)[2i],\Bbb Q(j)[2j])$ is 0 if $i\neq j$ and $\Bbb Q$ for $i=j$,\footnote{This follows since $M(\Bbb P^n )=\oplus_{i\in [0,n]} \Bbb Q(i)[2i]$ and $\End (M(\Bbb P^n))= CH_n (\Bbb P^n \times \Bbb P^n )=\Bbb Q^{[0,n]}$.  } and so every object  $M\in \CM (k)$ which is isomorphic to a direct sum of motives $\Bbb Q(i)[2i]$, $i\in \Bbb Z$, can be written in a unique manner as $\oplus_i \, V_i (i)[2i]$ where $V_i$ is a vector space (then $V_i =\Hom (\Bbb Q (i)[2i],M)$). Set $\tau^{\le 2a} M:= \oplus_{i\ge -a}\, V_i (i)[2i]$, etc.

\enspace

We are in the situation of  3.2 in the setting of $k$-varieties. As in loc.cit.,   $\CI^+ := \pi_* \Bbb Q (0)_Y \in D\CM (X_0 )$ (so  $\CI^{+\CH}$ is the corresponding Hodge sheaf from loc.cit.) Since $Y$ is smooth and $\pi$ is proper one has a natural self-duality $\CI^{+*}=\CI^+ (n-1)[2n-2]$.

The t-structure in $D\CM$ is not available, so we {\it define} the motivic intersection cohomology sheaf $\CI$ using a motivic version of decomposition (3.2.1):

\proclaim  {Proposition} There is a natural orthogonal direct sum decomposition in $D\CM (X_0 )$ $$\CI^+ =\CI  \oplus\, \oplus_\alpha i_{x_\alpha  *} \tau^{[2, 2n-4]}\Gamma (P_{\alpha})  \tag 4.2.1$$ whose Hodge realization is  (3.2.1) 
\endproclaim

\demo{Proof} It repeats the proof in 3.2 (minus its last paragraph). Namely, we first define a natural orthogonal decomposition  $$\Gamma (Z_\alpha )= H^{n-2}_{\text{prim}}(Z_\alpha )[2-n] \oplus \tau^{\le 2n-4} \Gamma (P_{\alpha}) \tag 4.2.2$$  in $D\CM (x_\alpha )=D\CM(k_{x_\alpha})$ whose Hodge realization is (3.2.2).\footnote{ So $H^{n-2}_{\text{prim}}(Z_\alpha )$ is a notation for a motive whose Hodge realization is the primitive cohomology of $Z_\alpha$; its definition does not involve any cohomology. To construct it explicitly, choose a $k$-point $z$ in $P_\alpha \smallsetminus Z_\alpha$. Let $\pi_z : Z_\alpha \to \Bbb P^{n-2}$ be the corresponding projection; this is a finite map of degree $d_\alpha$. Then  $H^{n-2}_{\text{prim}}(Z_\alpha )$ is the kernel of the projector $d_\alpha^{-1}\pi_z^t \pi_z$ acting on $M(Z_\alpha )(2-n)[4-2n]$.   } The construction in loc.cit.~uses only basic  six functors functoriality, so we can repeat it literally in the motivic setting. Then we proceed to define (4.2.1) as in loc.cit. \hfill$\square$
\enddemo

Set $\CB := \oplus_\alpha \, i_{x_\alpha *} H^{n-2}_{\text{prim}}(Z_\alpha )[1-n]\in D\CM(X_0 )$.
The self-dualities of $\Gamma (Z_\alpha )$ and of $\CI^+$, and the above orthogonal decompositions yield natural self-dualities $$\CB^* \iso \CB (n-2 )[2n-2], \quad \CI^* \iso  \CI (n-1)[2n-2].  \tag 4.2.3$$ 
 
 \proclaim {4.3. Lemma} (i) The adjunction  $\chi:  \Bbb Q (0)_{X_0} \to \pi_* \Bbb Q (0)_Y =\CI^+$ takes values in  $\CI \subset \CI^+$. \newline (ii) One has $\CC one (\chi : \Bbb Q (0)_{X_0} \to \CI )=\CB [1]$. 
\endproclaim

\demo{Proof}  (i) Follows since $\Hom (\Bbb Q (0)_{X_0}, i_{x_\alpha  *} \tau^{[2, 2n-4]}\Gamma (P_{\alpha}))=\Hom (\Bbb Q (0), \tau^{[2, 2n-4]}\Gamma (P_{\alpha}))$ $  =0$. \newline (ii) Since $\chi |_U =\id_{\Bbb Q (0)_U}$ the cone $\CC one (\chi )$ is supported on $\{ x_\alpha \}$. Now $i^*_{x_\alpha} \CC one (\chi )  = \CC one ( i^*_{x_\alpha}(\chi) )$  equals $H^{n-2}_{\text{prim}}(Z_\alpha )[2-n]$ by (4.2.2) and the construction of $\CI$, q.e.d. 
 \hfill$\square$
\enddemo

\remark{Remark} Since $\Ext^i (\Bbb Q(0)_{X_0}, \Bbb Q (0)^*_{X_0}(1-n)[2-2n])=\Ext^i (\Bbb Q (0),M(X_0 )(1-n)[2-2n])=CH_{n-1}(X_0 ,-i)$ we see that $\Ext^0 =Z_{n-1}(X_0 )$ and $Ext^{\neq 0}=0$, i.e., one has  $\Hom  (\Bbb Q(0)_{X_0}, \Bbb Q (0)^*_{X_0}(1-n)[2-2n]) =Z_{n-1}(X_0 )=Z_{n-1}(U)$.
\endremark

\remark{Example} One has    $\chi^* \chi =\epsilon$ where $\epsilon : \Bbb Q(0)_{X_0}\to\Bbb Q (0)^*_{X_0}(1-n)[2-2n]$ is the map that corresponds to the sum of irreducible components cycle (it is enough to check the assertion on $U$ where it is obvious).
\endremark

\subhead 4.4 \endsubhead We are in the situation of 3.3 in the setting of $k$-varieties.  Consider the   functor
  $ \psi^{\text{un}}_{tf}: D\CM (X\smallsetminus X_0 )\to D\CM (X_0 )$. There is a canonical morphism   $\iota: i_{X_0}^* \to \psi^{\text{un}}_{tf} v^*$ of functors on $D\CM (X)$ and its Verdier dual   $\iota^* :\psi^{\text{un}}_{tf} v^* \to i^!_{X_0}$.  Therefore we have a motivic sheaf $\CR :=\psi^{\text{un}}_{tf} \Bbb Q(0)_{X\smallsetminus X_0}$ equipped with a natural self-duality     $\CR^* \iso \CR (n-1)[2n-2] $ and   mutually dual maps $\Bbb Q (0)_{X_0} \buildrel{\iota}\over\to \CR \buildrel{\iota^*}\over\to \Bbb Q (0)^*_{X_0} (1-n)[2-2n]$ that are isomorphisms over $U$.

Let $\partial :\CB \to \Bbb Q(0)_{X_0}$ be the boundary map of the triangle from 4.3(ii). Set $\alpha_\CR := \iota\partial : \CB \to \CR$, and let $\beta_\CR $ be $\alpha_\CR^*$ combined with the self-duality identifications for $\CR$ and $\CB$, so we have
$$\CB  \buildrel{\alpha_\CR}\over\lra   \CR  \buildrel{\beta_\CR}\over\lra \CB (-1) . \tag 4.4.1   $$

\proclaim{Lemma-construction}  The composition $\beta_\CR \alpha_\CR$ is homotopic to zero. In fact, there is a canonical up to a homotopy  $\kappa_\CR$ such that $d(\kappa_\CR)=\beta_\CR \alpha_\CR$. 
\endproclaim

\demo{Proof}
By  Remark and Example in 4.3 one has  $\beta_\CR\alpha_\CR =\partial^* \iota^* \iota\partial = \partial^*\epsilon\partial =\partial^* \chi^* \chi\partial  =(\chi\partial)^* \chi\partial$. Notice that $\chi\partial$ is homotopic to 0; choose a homotopy $\lambda$, $d(\lambda ) =\chi\partial$.  Now set $\kappa_\CR :=\lambda^* \chi\partial $. 

 Independence  of $\kappa_\CR$ up to a homotopy from the choice of $\lambda$: if $\lambda'$ is another homotopy as above, i.e.,  $d(\lambda )=d(\lambda' )$, then $\kappa'_\CR =\lambda^{\prime *}\chi\partial = \kappa_\CR + (\lambda' -\lambda ) \chi\partial =\kappa_\CR + d((\lambda -\lambda') \lambda)$.  \hfill$\square$
\enddemo

\remark{Remark}  Our $\kappa_\CR$ is self-dual up to  homotopy: Indeed, one has $\kappa_\CR^* =(\chi\partial)^*\lambda =\kappa_\CR + d(\lambda^* \lambda )$. 
\endremark

\subhead 4.5 \endsubhead Below we use the notation from 2.1, 2.2. We have defined
$(\alpha_\CR ,\beta_\CR ,\kappa_\CR )\in D\CM (X_0 )^{(2)}$.  It yields the objects $E_\CR:= E ( \alpha_\CR ,\beta_\CR ,\kappa_\CR ) \in D\CM (X_0 )$ and
$(\alpha_\CI ,\beta_\CI ,\kappa_\CI )$ $:= \sigma (\alpha_\CR, \beta_\CR, \kappa_\CR )\in D\CM (X_0 )^{(2)}$. 
As follows from Remark in 4.4 and the definitions, the above three objects are naturally self-dual.

\proclaim{Proposition}
There is a  homotopy equivalence  
$\theta: \CI \iso E_\CR $  such that  the maps $\beta_\CI \theta: \CI \to \CB [1]$,   $\theta^{-1}\alpha_\CI: \CB (-1) [-1]\to \CI$ are a morphism of the triangle in 4.3(ii) and its dual. Our $\theta$ is unitary, i.e., $\theta^* =\theta^{-1}$.
\endproclaim

\demo{Proof} Recall that we have a natural homotopy equivalence 
 $(\lambda,\chi ): \CC one (\partial : \CB\to \Bbb Q(0)_{X_0}) \iso \CI$ (see 4.3(ii)), and $E_\CR$ is  the direct sum $\CB [1] \oplus \CR \oplus \CB (-1)[-1 ] $ with $(\alpha_\CR ,-\kappa_\CR ,\beta_\CR )$ added to the differential (see 2.1). Our $\theta$ is the composition $\CI \buildrel{\sim}\over{\leftarrow} \CC one (\partial ) \buildrel{\theta'}\over\to E_\CR$ where $\theta'$ is the next morphism: its restriction to
  $\CB [1] \subset \CC one (\partial )$ identifies it with the first summand in $E_\CR$, and its restriction to $\Bbb Q(0)_{X_0} \subset \CC one (\partial )$ is $(0, \iota , -\lambda^* \chi )$.

 One has $\theta^* \theta =\id_\CI$: we need to check that $\theta^{\prime*}\rho\theta'= (\lambda ,\chi )^* (\lambda ,\chi ): \CC one (\partial )\to \CC one (\partial^* )$ where $\rho : E_\CR \iso E_\CR^* (1-n)[2-2n]$ is the self-duality for $E_\CR$. As follows from Remark in 4.4, $\rho$ is the matrix  with the self-dualities for $\CR$ and $\CB$'s on the diagonal and the only non-zero off-diagonal entry being $\lambda^* \lambda : \CB \to \CB^* (1-n) [2-2n]$.    The rest  is an immediate calculation. 
 
 The assertion that $\beta_\CI \theta$ is the morphism of the triangle in 4.3(ii) means that $\beta_\CI \theta'$ is the projection $\CC one (\partial )\to \CB [1]$ which is evident from the construction. The assertion that $\alpha_\CI \theta$ is dual to $\beta_\CI \theta$ follows from the unitarity of $\theta$ once we know that $\theta$ is a homotopy equivalence. Let us check it.
 
 Our $\theta'$ is a morphism $\CC one (\CB \to \Bbb Q (0)_{X_0}) \to \CC one (\CB \to \CC one (\beta_\CR )[-1])$ compatible with the projections to $\CB$, and so it is enough to check that the map $(\iota , -\lambda^* \chi ) : \Bbb Q (0)_{X_0} \to \CC one (\beta_\CR )[-1]$ is a homotopy equivalence. Since $\iota$ is a homotopy equivalence on $ U$, it is enough to check our claim after applying $i^*_{x_\alpha}$. 
 
 The story of section 3.3 uses only the six functors formalism and basic facts from 3.1, so it remains literally true in the motivic setting.  Consider the canonical homotopy equivalence
 $a: i^*_{x_\alpha} \CR \iso \Gamma (P_\alpha \smallsetminus Z_\alpha )$ of  (3.3.1). By
 the Verdier dual assertion to the lemma in 3.3, $a$
  identifies $i^*_{x_\alpha} (\beta_\CR )$  with minus the residue map $r: \Gamma (P_\alpha \smallsetminus Z_\alpha ) \to 
 H^{n-2}_{\text{prim}}(Z_\alpha )(-1)[1-n] \subset \Gamma (Z_\alpha )(-1)[-1]$. By (4.2.2) we have a split exact triangle $\Bbb Q (0)\to \Gamma (P_\alpha \smallsetminus Z_\alpha ) \buildrel{r}\over\to 
 H^{n-2}_{\text{prim}}(Z_\alpha )(-1)[1-n]$, so
 $a$ identifies $i^*_{x_\alpha}\CC one (\beta_\CR )[-1])$ with $\Bbb Q (0) \subset \Gamma (P_\alpha \smallsetminus Z_\alpha ) $. It follows directly from the construction of $a$ that $a i^*_{x_\alpha}(\iota )$ coincides with the latter embedding, and we are done.
 \hfill$\square$
\enddemo 
  
\subhead 4.6 \endsubhead {\it Proof of the theorem in 1.9.} We have $(\alpha_\CI ,\beta_\CI ,\kappa_\CI ) \in D\CM(X_0 )^{(2)}$, hence $\Gamma (\alpha_\CI ,\beta_\CI ,\kappa_\CI )\in D\CM^{(2)}$.  For two Bloch cycles $A$, $B$ of classes $cl_A , cl_B \in \Hom ( \Bbb Q (0) ,$ $ H^{n-2}_{\text{prim}}(Z_\alpha ) (m))$ we have 
$( cl_A^* ,cl_{B*}) \Gamma (\alpha_\CI ,\beta_\CI ,\kappa_\CI )\in\fE^\CM (\Gamma (\CI )(m-1 )[1-n])=\fE^\CM (\Gamma (\CI^+ )(m-1)[1-n])=\fE^\CM ( M )$ where $M:=M(Y )(-m) [-1-2m ]$. By the construction
 the Hodge realization embedding $\fE^\CM ( M ) \hra \fE^\CH (M) = \fE^\CH (H^m (Y)) $ identifies it with $E^\psi_{A,B}$ from 1.6, and we are done. \hfill$\square$

\end
\Refs{}
\widestnumber\key{XXXXX}

\ref\key A1
\by J.~Ayoub
\book Les six op\'erations de Grothendieck et le formalisme des cycles \'evanescents dans le monde motivique (I)
\bookinfo Ast\'erisque 314
\publ SMF
\yr 2007
\endref


\ref\key A2
\by J.~Ayoub
\book Les six op\'erations de Grothendieck et le formalisme des cycles \'evanescents dans le monde motivique (II)
\bookinfo Ast\'erisque 315
\publ SMF
\yr 2007
\endref



\ref\key B
\by A.~Beilinson 
\paper Height pairing between algebraic cycles
\inbook K-theory, Arithmetic and Geometry, Yu.~I.~Manin (Ed.)
\bookinfo Lect.~Notes in Math.~1289
\publ Springer
\yr 1987
\endref

\ref\key Bl1
\by S.~Bloch 
\paper Height pairings for algebraic cycles
\jour Journal of Pure and Applied Algebra
\vol 34
\yr 1984
\pages 119--145
\endref

\ref\key Bl2
\by S.~Bloch 
\paper Cycles and biextensions
\jour Contemporary Mathematics
\vol 83
\yr 1989
\pages 19--30
\endref

\ref\key BlJS
\by S.~Bloch, R.~de Jong, E.~Can Sert\~oz
\paper Heights on curves and limits of Hodge structures
\jour arXiv:2206.01220
\yr 2022
\endref

\ref\key CD
\by D.-C.~Cisinski, F.~D\'eglise 
\book Triangulated categories of mixed motives
\publ Springer
\bookinfo Springer Monographs in Mathematics
\yr 2019
\endref

\ref\key G
\by S.~Gorchinskiy
\paper Notes on the biextension of Chow groups
\inbook Motives and algebraic cycles
\bookinfo Fields Institute Commun.
\vol 56
\pages 111--148
\publ Amer.~Math.~Soc.
\yr 2009
\endref

\ref\key Il
\by L.~Illusie 
\paper Sur la formule de Picard-Lefschetz
\inbook Algebraic geometry 2000, Azumino
\bookinfo Advanced Studies in Pure Math 
\vol 36
\pages 249--268
\publ Mathematical Society of Japan
\yr 2002
\endref

\ref\key Iv
\by F.~Ivorra 
\paper Perverse, Hodge and motivic realizations of \'etale motives
\jour Coompositio Mathematica
\vol 152
\issue 6
\yr 2016
\pages 1237--1285
\endref



\endRefs
